\documentclass[12pt]{article}
\usepackage{latexsym,amssymb,amsmath,enumerate,float,geometry,cite,tikz}
\usetikzlibrary{backgrounds,scopes,calc}
\tikzstyle{vertex}=[circle, draw, inner sep=0pt, minimum size = 4pt, fill]

\newtheorem{theorem}{Theorem}
\newtheorem{corollary}[theorem]{Corollary}

\newtheorem{lemma}[theorem]{Lemma}
\newtheorem{conjecture}[theorem]{Conjecture}

\newtheorem{claim}{Claim}

\usepackage{setspace}\onehalfspace

\begin{document}

\title{Relating Domination, Exponential Domination, \\ and Porous Exponential Domination}
\author{Michael A. Henning$^1$, Simon J\"{a}ger$^2$ and Dieter Rautenbach$^2$}
\date{}
\maketitle
\vspace{-10mm}
\begin{center}
{\small
$^1$ Department of Pure and Applied Mathematics, University of Johannesburg,\\
Auckland Park, 2006, South Africa, \texttt{mahenning@uj.ac.za}\\[3mm]
$^2$
Institute of Optimization and Operations Research, Ulm University,\\
Ulm, Germany, \texttt{simon.jaeger,dieter.rautenbach@uni-ulm.de}}
\end{center}

\begin{abstract}
The domination number $\gamma(G)$ of a graph $G$,
its exponential domination number $\gamma_e(G)$,
and
its porous exponential domination number $\gamma_e^*(G)$
satisfy $\gamma_e^*(G)\leq \gamma_e(G)\leq \gamma(G)$.
We contribute results about the gaps in these inequalities as well as the graphs for which some of the inequalities hold with equality.
Relaxing the natural integer linear program whose optimum value is $\gamma_e^*(G)$,
we are led to the definition of the fractional porous exponential domination number $\gamma_{e,f}^*(G)$ of a graph $G$.
For a subcubic tree $T$ of order $n$, we show $\gamma_{e,f}^*(T)=\frac{n+2}{6}$ and $\gamma_e(T)\leq 2\gamma_{e,f}^*(T)$.
We characterize the two classes of subcubic trees $T$ with $\gamma_e(T)=\gamma_{e,f}^*(T)$ and $\gamma(T)=\gamma_e(T)$, respectively.
Using linear programming arguments, we establish several lower bounds on the fractional porous exponential domination number in more general settings.
\end{abstract}

\pagebreak

\section{Introduction}

In \cite{ddems} Dankelmann et al.~introduce {\it exponential domination} as a variant of domination in graphs
where the influence of the vertices in an {\it exponential dominating set}
extends to any arbitrary distance but decays exponentially with that distance.
They consider two parameters for a given graph $G$,
its {\it exponential domination number} $\gamma_e(G)$,
corresponding to a setting in which the different vertices in the exponential dominating set block each others influence, and
its {\it porous exponential domination number} $\gamma_e^*(G)$,
where such a blocking does not occur.

Unlike most other domination parameters \cite{hhs}, which are based on local conditions,
exponential domination is a genuinely global concept.
Compared to exponential domination,
even notions such as distance domination \cite{he} appear as essentially local,
because they can be reduced to ordinary domination by considering suitable powers of the underlying graph.
The global nature of exponential domination makes it much harder,
which might be the reason why there are only relatively few results about it \cite{abcvy,aa,bor1,bor2}.
While Bessy et al.~\cite{bor2} show that the exponential domination number is APX-hard for subcubic graphs
and describe an efficient algorithm for subcubic trees,
the complexity of the exponential domination number on general trees is unknown,
and a hardness result does not seem unlikely.
For the porous version, even less is known.
There is not a single complexity result,
and it is even unknown whether the porous exponential domination number of a given subcubic tree can be determined efficiently.
Partly motivated by such difficulties,
Goddard et al.~define \cite{ghp} {\it disjunctive domination} (see also \cite{hm,hm2,ppp}),
which keeps the exponential decay of the influence but
only considers distances one and two.
Further related parameters known as {\it influence} and {\it total influence} \cite{dll} also have unknown complexity even for trees \cite{hed}.

The two parameters of exponential domination and the classical domination number $\gamma(G)$ \cite{hhs} of a graph $G$ satisfy
$$\gamma_e^*(G)\leq \gamma_e(G)\leq \gamma(G).$$
In the present paper we contribute results about the gaps in these inequalities
as well as
the graphs for which some of the inequalities hold with equality.
In order to obtain lower bounds
we consider fractional relaxations and apply linear programming techniques.

\medskip

\noindent Before stating our results and several open problems, we recall some notation and give formal definitions.
We consider finite, simple, and undirected graphs.
The vertex set and the edges set of a graph $G$ are denoted by $V(G)$ and $E(G)$, respectively.
The order $n(G)$ of $G$ is the number of vertices of $G$.
If all vertex degrees in $G$ are at most $3$, then $G$ is subcubic.
A vertex of degree at most $1$ in $G$ is an endvertex of $G$.
A vertex in a rooted tree $T$ is a leaf of $T$ if it has no child in $T$.
The distance ${\rm dist}_G(u,v)$ between two vertices $u$ and $v$ in $G$
is the minimum number of edges of a path in $G$ between $u$ and $v$.
If no such path exists, then let ${\rm dist}_G(u,v)=\infty$.
The diameter ${\rm diam}(G)$ of $G$ is the maximum distance between vertices of $G$.
A set $D$ of vertices of a graph $G$ is a dominating set of $G$ \cite{hhs}
if every vertex of $G$ not in $D$ has a neighbor in $D$,
and the domination number $\gamma(G)$ of $G$ is the minimum order of a dominating set of $G$.
Similarly, for some set $X$ of vertices of $G$,
let $\gamma(G,X)$ be the minimum order of a set $D$ of vertices such that every vertex in $X\setminus D$ has a neighbor in $D$.
Note that $\gamma(G)=\gamma(G,V(G))$.
For positive integers $n$ and $m$,
let $[n]$ be the set of the positive integers at most $n$,
let $P_n$, $C_n$, and $K_n$ be the path, cycle, and complete graph of order $n$, respectively,
and, let $K_{n,m}$ be the complete bipartite graph with partite sets of orders $n$ and $m$.

In order to capture the above-mentioned blocking effects that are a feature of exponential domination,
we need a modified distance notion.
Therefore, let $D$ be a set of vertices of a graph $G$.
For two vertices $u$ and $v$ of $G$,
let ${\rm dist}_{(G,D)}(u,v)$ be the minimum number of edges of a path $P$ in $G$ between $u$ and $v$
such that $D$ contains exactly one endvertex of $P$ but no internal vertex of $P$.
If no such path exists, then let ${\rm dist}_{(G,D)}(u,v)=\infty$.
Note that, if $u$ and $v$ are distinct vertices in $D$,
then ${\rm dist}_{(G,D)}(u,u)=0$ and ${\rm dist}_{(G,D)}(u,v)=\infty$.

For a vertex $u$ of $G$, let
\begin{eqnarray}\label{ew}
w_{(G,D)}(u)=\sum\limits_{v\in D}\left(\frac{1}{2}\right)^{{\rm dist}_{(G,D)}(u,v)-1},
\end{eqnarray}
where $\left(\frac{1}{2}\right)^{\infty}=0$.
Note that $w_{(G,D)}(u)=2$ for $u\in D$.

Dankelmann et al.~\cite{ddems} define
the set $D$ to be an {\it exponential dominating set} of $G$ if
$$\mbox{$w_{(G,D)}(u)\geq 1$ for every vertex $u$ of $G$,}$$
and the {\it exponential domination number $\gamma_e(G)$} of $G$
as the minimum order of an exponential dominating set.
Similarly, they define $D$ to be a {\it porous exponential dominating set} of $G$ if
$$\mbox{$w^*_{(G,D)}(u)\geq 1$ for every vertex $u$ of $G$,}$$
where
\begin{eqnarray}\label{ew*}
w^*_{(G,D)}(u)=\sum\limits_{v\in D}\left(\frac{1}{2}\right)^{{\rm dist}_G(u,v)-1},
\end{eqnarray}
and they define the {\it porous exponential domination number $\gamma^*_e(G)$} of $G$
as the minimum order of a porous exponential dominating set of $G$.
Note that the definition of $w^*_{(G,D)}(u)$ involves the usual distance rather than ${\rm dist}_{(G,D)}(u,v)$,
which reflects the absence of blocking effects.
A dominating, exponential dominating, or porous exponential dominating set of minimum order is called {\it minimum}.

\medskip

\noindent The parameter $\gamma^*_e(G)$ equals the optimum value of the following integer linear program
\begin{eqnarray}\label{epi}
\begin{array}{lrlll}
\min	& \sum\limits_{u\in V(G)}x(u) & & & \\[5mm]
s.t.		& \sum\limits_{u\in V(G)}\left(\frac{1}{2}\right)^{{\rm dist}_{G}(u,v)-1}\cdot x(u) & \geq & 1 & \forall v\in V(G)\\
		& x(u) & \in & \{ 0,1\} & \forall u\in V(G).
\end{array}
\end{eqnarray}
Relaxing (\ref{epi}), we obtain the following linear program
\begin{eqnarray}\label{ep}
\begin{array}{lrlll}
\min	& \sum\limits_{u\in V(G)}x(u) & & & \\[5mm]
s.t.		& \sum\limits_{u\in V(G)}\left(\frac{1}{2}\right)^{{\rm dist}_{G}(u,v)-1}\cdot x(u) & \geq & 1 & \forall v\in V(G)\\
		& x(u) & \geq & 0 & \forall u\in V(G).
\end{array}
\end{eqnarray}
Let the {\it fractional porous exponential domination number $\gamma^*_{e,f}(G)$} of $G$ be the optimum value of (\ref{ep}).

Clearly,
\begin{eqnarray}\label{eseq}
\gamma_{e,f}^*(G)\leq \gamma_e^*(G)\leq \gamma_e(G)\leq \gamma(G)
\end{eqnarray}
for every graph $G$.

Most of our results concern subcubic graphs.
While exponential domination number is APX-hard already for subcubic graphs~\cite{bor2}, 
there is the following fundamental lemma from \cite{bor2},
which will be an important technical tool throughout our paper.

\begin{lemma}[Bessy et al.~\cite{bor2}]\label{lemma1}
Let $G$ be a graph of maximum degree at most $3$,
and let $D$ be a set of vertices of $G$.

If $u$ is a vertex of degree at most $2$ in $G$, then $w_{(G,D)}(u)\leq 2$ with equality if and only if
$u$ is contained in a subgraph $T$ of $G$ that is a tree,
such that rooting $T$ in $u$ yields a full binary tree and $D\cap V(T)$ is exactly the set of leaves of $T$.
\end{lemma}
The next section contains all our results as well as many closely related conjectures and open problems.

\section{Results}

Our first slightly surprising result is that the fractional porous exponential domination number of a subcubic tree only depends on its order.

\begin{theorem}\label{theorem2}
If $T$ is a subcubic tree of order $n$, then $\gamma_{e,f}^*(T)=\frac{n+2}{6}$.
\end{theorem}
{\it Proof:}
Let $T$ be a subcubic tree of order $n$.
If $T$ has only one vertex $u$, then $x(u)=\frac{1}{2}=\frac{1+2}{6}$ is an optimum solution of (\ref{ep}).
Hence, we may assume that $n\geq 2$.
Let $V_i$ be the set of vertices of degree $i$ in $T$, and let $n_i=|V_i|$ for $i\in [3]$.
Let $T'$ arise from $T$ by adding, for every vertex $u$ in $V_2$,
a new vertex $p_u$ as well as the new edge $up_u$.
By construction, $T'$ is a tree of order $n+n_2$ that only has vertices of degree $1$ and $3$,
and $D'=V_1\cup \{ p_u:u\in V_2\}$ is the set of all endvertices of $T'$.

Let $v$ be an endvertex of $T'$, and let $u$ be the neighbor of $v$.
Since $T'-v$ rooted in $u$ is a full binary tree whose set of leaves is exactly $D'\setminus \{ v\}$,
Lemma \ref{lemma1} implies that $w_{(T'-v,D'\setminus \{ v\})}(u)=2$.
Since $v\in D'$, this implies $w_{(T',D')}(v)=2+1=3$.
Similarly, if $v$ is a vertex of degree $3$ in $T'$ whose neighbors are $u_1$, $u_2$, and $u_3$,
then Lemma \ref{lemma1} implies that $w_{(T'-v,D')}(u_i)=2$ for $i\in [3]$,
which implies $w_{(T',D')}(v)=1+1+1=3$.
Since $D'$ only contains leaves of $T'$,
we obtain that $w_{(T',D')}^*(v)=w_{(T',D')}(v)=3$ holds for every vertex $v$ of $T'$.

Let $(x(u))_{u\in V(T)}$ be such that
$$x(u)
=
\begin{cases}
\frac{1}{3} & \mbox{, if $u$ is an endvertex of $T$,}\\
\frac{1}{6} & \mbox{, if $u$ has degree $2$ in $T$, and}\\
0 & \mbox{, if $u$ has degree $3$ in $T$}.
\end{cases}
$$
If $u$ and $v$ are vertices of $T$, then ${\rm dist}_{T}(u,v)={\rm dist}_{T'}(u,v)$.
Furthermore, if $u\in V_2$, then ${\rm dist}_{T}(u,v)={\rm dist}_{T'}(p_u,v)-1$.

This implies that
\begin{eqnarray}
\sum\limits_{u\in V(T)}\left(\frac{1}{2}\right)^{{\rm dist}_{T}(u,v)-1}\cdot x(u)
& = &
\frac{1}{3}\sum\limits_{u\in V_1}\left(\frac{1}{2}\right)^{{\rm dist}_{T}(u,v)-1}
+\frac{1}{6}\sum\limits_{u\in V_2}\left(\frac{1}{2}\right)^{{\rm dist}_{T}(u,v)-1}\nonumber\\
& = &
\frac{1}{3}\sum\limits_{u\in V_1}\left(\frac{1}{2}\right)^{{\rm dist}_{T'}(u,v)-1}
+\frac{1}{3}\sum\limits_{u\in V_2}\left(\frac{1}{2}\right)^{{\rm dist}_{T'}(p_u,v)-1}\nonumber\\
& = &
\frac{1}{3}\sum\limits_{u\in D'}\left(\frac{1}{2}\right)^{{\rm dist}_{T'}(u,v)-1}\nonumber\\
& = & \frac{1}{3}w^*_{(T',D')}(v)\nonumber\\
& = & 1\label{eref1}
\end{eqnarray}
for every vertex $v$ of $T$.

The dual linear program of (\ref{ep}) is
\begin{eqnarray}\label{ed}
\begin{array}{lrlll}
\max	& \sum\limits_{v\in V(T)}y(v) & & & \\[5mm]
s.t.		& \sum\limits_{v\in V(T)}\left(\frac{1}{2}\right)^{{\rm dist}_{T}(u,v)-1}\cdot y(v) & \leq & 1 & \forall u\in V(T)\\
		& y(v) & \geq & 0 & \forall v\in V(T).
\end{array}
\end{eqnarray}
Therefore, setting $y(u)=x(u)$ for every vertex $u$ of $T$,
we obtain, by (\ref{eref1}), that
\begin{itemize}
\item $(x(u))_{u\in V(T)}$ is a feasible solution of (\ref{ep}),
\item $(y(u))_{u\in V(T)}$ is a feasible solution of (\ref{ed}), and that
\item $\sum\limits_{u\in V(T)}x(u)=\sum\limits_{u\in V(T)}y(u)$,
\end{itemize}
that is, $(x(u))_{u\in V(T)}$ and $(y(u))_{u\in V(T)}$ are both optimum solutions of the respective linear programs.

Since $n_1=n_3+2$, we obtain $n+2=n_1+n_2+n_3+2=2n_1+n_2$ and, hence,
$$\gamma_{e,f}^*(T)=\sum\limits_{u\in V(T)}x(u)
=\frac{1}{3}n_1+\frac{1}{6}n_2
=\frac{n+2}{6},$$
which completes the proof.
$\Box$

\medskip

\noindent Bessy et al.~\cite{bor2} show
$\frac{n+2}{6}\leq \gamma_e(T)\leq \frac{n+2}{3}$
for every subcubic tree $T$ of order $n$.
Note that the first of these two inequalities is an immediate consequence of Theorem \ref{theorem2},
and that,
combined with Theorem \ref{theorem2},
the second of these inequalities implies the following.

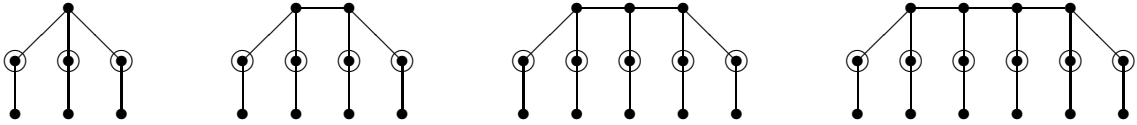
\begin{figure}[H]
\begin{center}
$\mbox{}$\hfill
\unitlength 0.7mm 
\linethickness{0.4pt}
\ifx\plotpoint\undefined\newsavebox{\plotpoint}\fi 
\begin{picture}(27,26)(0,0)
\put(5,5){\circle*{2}}
\put(15,5){\circle*{2}}
\put(25,5){\circle*{2}}
\put(5,15){\circle*{2}}
\put(15,15){\circle*{2}}
\put(25,15){\circle*{2}}
\put(15,25){\circle*{2}}
\put(25,5){\line(0,1){10}}
\put(25,15){\line(-1,1){10}}
\put(15,25){\line(0,-1){20}}
\put(15,25){\line(-1,-1){10}}
\put(5,15){\line(0,-1){10}}
\put(5,15){\circle{4}}
\put(15,15){\circle{4}}
\put(25,15){\circle{4}}
\end{picture}
\hfill
\linethickness{0.4pt}
\ifx\plotpoint\undefined\newsavebox{\plotpoint}\fi 
\begin{picture}(37,26)(0,0)
\put(5,5){\circle*{2}}
\put(25,5){\circle*{2}}
\put(15,5){\circle*{2}}
\put(35,5){\circle*{2}}
\put(5,15){\circle*{2}}
\put(25,15){\circle*{2}}
\put(15,15){\circle*{2}}
\put(35,15){\circle*{2}}
\put(15,25){\circle*{2}}
\put(15,25){\line(0,-1){20}}
\put(15,25){\line(-1,-1){10}}
\put(5,15){\line(0,-1){10}}
\put(25,15){\line(0,-1){10}}
\put(5,15){\circle{4}}
\put(25,15){\circle{4}}
\put(15,15){\circle{4}}
\put(35,15){\circle{4}}
\put(35,15){\line(0,-1){10}}
\put(25,25){\circle*{2}}
\put(25,24){\line(0,-1){9}}
\put(25,25){\line(1,-1){10}}
\put(15,25){\line(1,0){9}}
\end{picture}
\hfill
\linethickness{0.4pt}
\ifx\plotpoint\undefined\newsavebox{\plotpoint}\fi 
\begin{picture}(47,26)(0,0)
\put(5,5){\circle*{2}}
\put(35,5){\circle*{2}}
\put(15,5){\circle*{2}}
\put(25,5){\circle*{2}}
\put(45,5){\circle*{2}}
\put(5,15){\circle*{2}}
\put(35,15){\circle*{2}}
\put(15,15){\circle*{2}}
\put(25,15){\circle*{2}}
\put(45,15){\circle*{2}}
\put(15,25){\circle*{2}}
\put(25,25){\circle*{2}}
\put(15,25){\line(0,-1){20}}
\put(25,25){\line(0,-1){20}}
\put(15,25){\line(-1,-1){10}}
\put(5,15){\line(0,-1){10}}
\put(35,15){\line(0,-1){10}}
\put(5,15){\circle{4}}
\put(35,15){\circle{4}}
\put(15,15){\circle{4}}
\put(25,15){\circle{4}}
\put(45,15){\circle{4}}
\put(45,15){\line(0,-1){10}}
\put(35,25){\circle*{2}}
\put(35,24){\line(0,-1){9}}
\put(35,25){\line(1,-1){10}}
\put(15,25){\line(1,0){9}}
\put(25,25){\line(1,0){9}}
\end{picture}
\hfill
\linethickness{0.4pt}
\ifx\plotpoint\undefined\newsavebox{\plotpoint}\fi 
\begin{picture}(57,26)(0,0)
\put(5,5){\circle*{2}}
\put(45,5){\circle*{2}}
\put(15,5){\circle*{2}}
\put(25,5){\circle*{2}}
\put(35,5){\circle*{2}}
\put(55,5){\circle*{2}}
\put(5,15){\circle*{2}}
\put(45,15){\circle*{2}}
\put(15,15){\circle*{2}}
\put(25,15){\circle*{2}}
\put(35,15){\circle*{2}}
\put(55,15){\circle*{2}}
\put(15,25){\circle*{2}}
\put(25,25){\circle*{2}}
\put(35,25){\circle*{2}}
\put(15,25){\line(0,-1){20}}
\put(25,25){\line(0,-1){20}}
\put(35,25){\line(0,-1){20}}
\put(15,25){\line(-1,-1){10}}
\put(5,15){\line(0,-1){10}}
\put(45,15){\line(0,-1){10}}
\put(5,15){\circle{4}}
\put(45,15){\circle{4}}
\put(15,15){\circle{4}}
\put(25,15){\circle{4}}
\put(35,15){\circle{4}}
\put(55,15){\circle{4}}
\put(55,15){\line(0,-1){10}}
\put(45,25){\circle*{2}}
\put(45,24){\line(0,-1){9}}
\put(45,25){\line(1,-1){10}}
\put(15,25){\line(1,0){9}}
\put(25,25){\line(1,0){9}}
\put(35,25){\line(1,0){9}}
\end{picture}
\hfill$\mbox{}$
\end{center}
\caption{Lemma \ref{lemma1} easily implies that every exponential dominating set of the trees shown above intersects the closed neighborhood of each endvertex.
Therefore, the circled vertices form a minimum exponential dominating set,
and hence $\gamma_e(T)=\frac{n(T)+2}{3}=2\gamma_{e,f}^*(T)$ for each shown tree $T$.
Furthermore, it is easy to verify that
$\gamma_e^*(T)=\gamma_e(T)$
for the two smallest trees $T$ of orders $7$ and $10$.}\label{fign1}
\end{figure}

\begin{corollary}\label{corollary1}
If $T$ is a subcubic tree of order $n$, then $\gamma_e(T)\leq 2\gamma_{e,f}^*(T)$.
\end{corollary}
Figure \ref{fign1} illustrates an infinite family of trees showing that Corollary \ref{corollary1} is tight.
Another immediate consequence of Corollary \ref{corollary1},
namely $\gamma_e^*(T)\leq 2\gamma_{e,f}^*(T)$ for every subcubic tree $T$,
is tight for the two smallest trees in this family,
which implies that the integrality gap between the integer linear program (\ref{epi}) and its linear programming relaxation (\ref{ep}) is $2$ for such trees.
It seems possible that the integrality gap between
(\ref{epi}) and (\ref{ep}) is bounded for all graphs of bounded maximum degree.

\noindent Yet another consequence of Corollary \ref{corollary1} is that $\gamma_e(T)\leq 2\gamma_e^*(T)$
for every subcubic tree $T$.
We believe that this estimate can be improved as follows.

\begin{figure}[H]
\begin{center}
\unitlength 0.7mm 
\linethickness{0.4pt}
\ifx\plotpoint\undefined\newsavebox{\plotpoint}\fi 
\begin{picture}(116,42)(0,0)
\put(5,5){\circle*{2}}
\put(45,5){\circle*{2}}
\put(85,5){\circle*{2}}
\put(25,5){\circle*{2}}
\put(65,5){\circle*{2}}
\put(105,5){\circle*{2}}
\put(15,5){\circle*{2}}
\put(55,5){\circle*{2}}
\put(95,5){\circle*{2}}
\put(35,5){\circle*{2}}
\put(75,5){\circle*{2}}
\put(115,5){\circle*{2}}
\put(10,15){\circle*{2}}
\put(50,15){\circle*{2}}
\put(90,15){\circle*{2}}
\put(30,15){\circle*{2}}
\put(70,15){\circle*{2}}
\put(110,15){\circle*{2}}
\put(20,25){\circle*{2}}
\put(60,25){\circle*{2}}
\put(100,25){\circle*{2}}
\put(5,5){\line(1,2){5}}
\put(45,5){\line(1,2){5}}
\put(85,5){\line(1,2){5}}
\put(10,15){\line(1,1){10}}
\put(50,15){\line(1,1){10}}
\put(90,15){\line(1,1){10}}
\put(20,25){\line(1,-1){10}}
\put(60,25){\line(1,-1){10}}
\put(100,25){\line(1,-1){10}}
\put(30,15){\line(1,-2){5}}
\put(70,15){\line(1,-2){5}}
\put(110,15){\line(1,-2){5}}
\put(30,15){\line(-1,-2){5}}
\put(70,15){\line(-1,-2){5}}
\put(110,15){\line(-1,-2){5}}
\put(10,15){\line(1,-2){5}}
\put(50,15){\line(1,-2){5}}
\put(90,15){\line(1,-2){5}}
\put(60,40){\circle*{2}}
\multiput(100,25)(-.0898876404,.0337078652){445}{\line(-1,0){.0898876404}}
\put(60,40){\line(0,-1){15}}
\multiput(60,40)(-.0898876404,-.0337078652){445}{\line(-1,0){.0898876404}}
\put(10,15){\circle{4}}
\put(30,15){\circle{4}}
\put(50,15){\circle{4}}
\put(70,15){\circle{4}}
\put(90,15){\circle{4}}
\put(110,15){\circle{4}}
\put(18,23){\framebox(4,4)[cc]{}}
\put(58,23){\framebox(4,4)[cc]{}}
\put(58,38){\framebox(4,4)[cc]{}}
\put(98,23){\framebox(4,4)[cc]{}}
\end{picture}
\end{center}
\caption{A tree $T$ with $\gamma_e(T)=6$ and $\gamma_e^*(T)=\gamma_{e,f}^*(T)=4$.
The circled vertices indicate a minimum exponential dominating set
while the boxed vertices indicate a minimum porous exponential dominating set.}\label{figure1}
\end{figure}
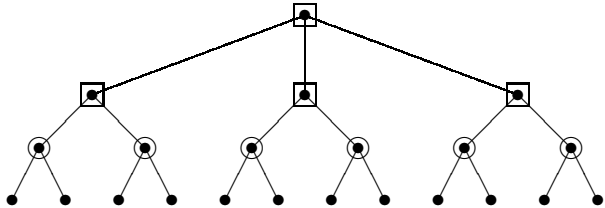

\begin{conjecture}\label{conjecture1}
If $T$ is a subcubic tree, then $\gamma_e(T)\leq \frac{3}{2}\gamma^*_e(T)$.
\end{conjecture}
The tree in Figure \ref{figure1} shows that the bound in Conjecture \ref{conjecture1} would be tight.

\noindent The special tree in Figure \ref{figure1}
also plays a role in our next conjecture.

\begin{conjecture}\label{conjecture2}
If $T$ is a subcubic tree, then $\gamma_e^*(T)=\gamma_{e,f}^*(T)$
if and only if $T$ is either $K_{1,3}$ or the tree in Figure \ref{figure1}.
\end{conjecture}
We establish many quite restrictive properties of the trees considered in Conjecture \ref{conjecture2}.

\begin{theorem}\label{theorem3}
Let $T$ be a subcubic tree with $\gamma_e^*(T)=\gamma_{e,f}^*(T)>1$.
Let $V_i$ be the set of vertices of degree $i$ in $T$ for $i\in [3]$. 
Let $D$ be a minimum porous exponential dominating set of $T$.
\begin{enumerate}[(i)]
\item $w_{(T,D)}^*(u)=1$ for every vertex $u\in V_1\cup V_2$.
\item $D\subseteq V_3$ and $N_T(V_1\cup V_2)\subseteq V_3\setminus D$.
\item $T$ does not contain a vertex $u$ in $V_1$ and a vertex $w$ in $V_2$ at distance $2$.
\item $T$ does not contain two vertices $u_1$ and $u_2$ in $V_1$ and a vertex $v$ in $V_2$
such that ${\rm dist}_T(u_1,u_2)=2$ and ${\rm dist}_T(u_1,v)\in \{ 3,4\}$.
\end{enumerate}
\end{theorem}
{\it Proof:} Note that $\gamma_e^*(T)>1$ implies that $T$ is not a star.

\medskip

\noindent (i) Let $(x(u))_{u\in V(T)}$ be such that
$$x(u)
=
\begin{cases}
1 & \mbox{, if $u\in D$, and }\\
0 & \mbox{, otherwise}
\end{cases}
$$
and let $(y(u))_{u\in V(T)}$ be as in the proof of Theorem \ref{theorem2}.
Since $\gamma_e^*(T)=\gamma_{e,f}^*(T)$,
we obtain that
$(x(u))_{u\in V(T)}$ is an optimum solution of (\ref{ep}), and
$(y(u))_{u\in V(T)}$ is an optimum solution of (\ref{ed}).
By the dual complementary slackness conditions,
we obtain that $y(v)>0$ for some $v\in V(T)$ implies
$$w_{(T,D)}^*(v)=\sum\limits_{u\in V(G)}\left(\frac{1}{2}\right)^{{\rm dist}_{T}(u,v)-1}\cdot x(u)=1.$$
Since $y(v)>0$ if and only if $v\in V_1\cup V_2$,
the proof of (i) is complete.

\medskip

\noindent (ii) Since $v\in D$ implies $w_{(T,D)}^*(v)\geq 2$, (i) implies $D\subseteq V_3$.
If $u$ in $V_1$ has a neighbor $v$ in $V_2$, then $w_{(T,D)}^*(v)=1$ implies the contradiction $w_{(T,D)}^*(u)=2$.
Hence, $N_T(V_1)\subseteq V_3$.

Suppose that $T$ contains a path $uvwx$ with $d_T(v)=d_T(w)=2$.
Let $T_v$ be the component of $T-w$ that contain $v$.
Let $\alpha=w^*_{(T,D\cap V(T_v))}(v)$, and $\beta=w^*_{(T,D\setminus V(T_v))}(v)$.
Clearly,
$w^*_{(T,D\cap V(T_v))}(w)=\frac{1}{2}\alpha$, and
$w^*_{(T,D\setminus V(T_v))}(w)=2\beta$.
Since $w_{(T,D)}^*(v)=\alpha+\beta=1$ and $w_{(T,D)}^*(w)=\frac{1}{2}\alpha+2\beta=1$,
we obtain $\beta=\frac{1}{3}$.
By the definition (\ref{ew*}) of $w^*$, we obtain that $\beta$ is the finite sum of powers of $2$,
that is, $\beta=2^{p_1}+\ldots+2^{p_k}$ for suitable integers $p_1,\ldots,p_k$.
Let $p=\max\{ |p_1|,\ldots,|p_k|\}$.
Since $\frac{1}{3}=2^{p_1}+\ldots+2^{p_k}$,
we obtain $2^p=3\cdot (2^{p_1}+\ldots+2^{p_k})\cdot 2^p$.
Note that $(2^{p_1}+\ldots+2^{p_k})\cdot 2^p$ is an integer.
While the right hand side of this equation is divisible by $3$,
the left hand side is not, which is a contradiction.
Hence, $N_T(V_2)\subseteq V_3$.

If $u\in V_1\cup V_2$ has a neighbor in $D$, then $|D|=\gamma_e^*(T)>1$ implies the contradiction $w_{(T,D)}^*(u)>1$.
Hence, $N_T(V_1\cup V_2)\subseteq V_3\setminus D$.

\medskip

\noindent (iii) Let $v$ be the common neighbor of $u$ and $w$.
By (ii), we have $v\in V_3\setminus D$.
Let $T_v$ be the component of $T-w$ that contains $v$.
Let $\alpha=w_{(T,D\cap V(T_v))}^*(v)$ and $\beta=w_{(T,D\setminus V(T_v))}^*(w)$.
By (i), we obtain
$w_{(T,D)}^*(u)=\frac{1}{2}\alpha+\frac{1}{4}\beta=1$
and
$w_{(T,D)}^*(w)=\frac{1}{2}\alpha+\beta=1$,
which implies $\alpha=2$ and $\beta=0$.
Since $\beta=0$, we obtain $D\setminus V(T_v)=\emptyset$.
Now, if $u'$ is an endvertex of $T$ in the component of $T-v$ that contains $w$, then $w_{(T,D)}^*(u')\leq \frac{1}{4}\alpha<1$,
which is a contradiction.

\medskip

\noindent (iv)
First, we assume that ${\rm dist}_T(u_1,v)=3$.
Let $u_1xyv$ be a path in $T$.
By (ii), the vertices $x$ and $y$ belong to $V_3\setminus D$.
Let $T_y$ be the component of $T-v$ that contains $y$.
Let
$\alpha=w_{(T,D\cap V(T_y))}^*(y)$ and
$\beta=w_{(T,D\setminus V(T_y))}^*(v)$.
By (i), we obtain
$w_{(T,D)}^*(u_1)=w_{(T,D)}^*(u_2)=\frac{1}{4}\alpha+\frac{1}{8}\beta=1$
and
$w_{(T,D)}^*(v)=\frac{1}{2}\alpha+\beta=1$,
which implies
$\frac{1}{4}\alpha+\frac{7}{8}\beta
=\left(\frac{1}{2}\alpha+\beta\right)-\left(\frac{1}{4}\alpha+\frac{1}{8}\beta\right)
=0$.
Since, by definition, $\alpha$ and $\beta$ are non-negative,
we obtain the contradiction $\alpha=\beta=0$.

Next, we assume that ${\rm dist}_T(u_1,v)=4$.
Let $u_1xyzv$ be a path in $T$.
Let $T_y$ be the component of $T-\{ xy,yz\}$ that contains $y$, and
let $T_z$ be the component of $T-\{ yz,zv\}$ that contains $z$.
Let
$\alpha=w_{(T,D\cap V(T_y))}^*(y)$,
$\beta=w_{(T,D\cap V(T_z))}^*(z)$, and
$\gamma=w_{(T,D\setminus (V(T_y)\cup V(T_z)))}^*(v)$.
By (i), we obtain
$w_{(T,D)}^*(u_1)=w_{(T,D)}^*(u_2)=\frac{1}{4}\alpha+\frac{1}{8}\beta+\frac{1}{16}\gamma=1$
and
$w_{(T,D)}^*(v)=\frac{1}{4}\alpha+\frac{1}{2}\beta+\gamma=1$,
which implies
$\frac{3}{8}\beta+\frac{15}{16}\gamma=0$,
and, hence, $\beta=\gamma=0$, and $\alpha=4$.
Now, if $u$ is an endvertex of $T$ that lies in the component of $T-z$ that contains $v$,
then $w_{(T,D)}^*(u)\leq \frac{1}{8}\alpha<1$,
which is a contradiction.
$\Box$

\medskip

\noindent Our next result would be an immediate consequence of Conjecture \ref{conjecture2}.

\begin{theorem}\label{theorem4}
If $T$ is a subcubic tree, then $\gamma_e(T)=\gamma_{e,f}^*(T)$
if and only if $T$ is $K_{1,3}$.
\end{theorem}
{\it Proof:} If $T=K_{1,3}$, then $\gamma_e(T)=1=\frac{4+2}{6}=\gamma_{e,f}^*(T)$,
which implies the sufficiency.
In order to prove the necessity, let $T$ be a subcubic tree with $\gamma_e(T)=\gamma_{e,f}^*(T)$.
By~(\ref{eseq}), $\gamma_e^*(T) = \gamma_e(T)$.
If $\gamma_e(T)=1$, then $T\in \{ P_1,P_2,P_3,K_{1,3}\}$.
By Theorem \ref{theorem2}, the only tree $T$ in this set with $\gamma_{e,f}^*(T)=1$ is $K_{1,3}$.
Hence, we may assume that $\gamma_e(T)>1$, which implies $\gamma^*_e(T)>1$.
Let $D$ be a minimum exponential dominating set.
Let $u$ be an endvertex of $T$.
Since $D$ is also a minimum porous exponential dominating set,
Theorem \ref{theorem3} implies $w_{(T,D)}^*(u)=1$.
Now, $1\leq w_{(T,D)}(u)\leq w_{(T,D)}^*(u)=1$, which implies $w_{(T,D)}(u)=w_{(T,D)}^*(u)=1$.
If $v$ is the neighor of $u$, then $w_{(T,D)}(v)=2$.
Since $w_{(T,D)}(u)=w_{(T,D)}^*(u)$ and $|D|\geq 2$, this implies $v\not\in D$.
By Lemma \ref{lemma1},
the tree $T$ contains a full binary tree $T'$ rooted in $v$
such that $V(T')\cap D$ is exactly the set of leaves of $T'$.
Since $u\not\in D$, the vertex $u$ does not belong to $T'$.
Since $v\not\in D$, the tree $T'$ has at least two leaves. 
Let $u'$ be a leaf of $T'$, and let $v'$ be the parent of $u'$ in $T'$.
By Theorem \ref{theorem3}(ii), 
$D \subseteq V_3$, implying that no endvertex of $T'$ is a leaf of $T$. In particular, $u'$ is not an endvertex of $T'$.
Let $u''$ be an endvertex of $T$ in the component of $T-v'$ that contains $u'$.
Since $T'$ has more than one leaf, we obtain $w_{(T,D)}(u'')<w_{(T,D)}^*(u'')$.
Nevertheless, arguing as above we obtain $w_{(T,D)}^*(u'')=1$,
which implies the contradiction
$1\leq w_{(T,D)}(u'')<w_{(T,D)}^*(u'')=1$.
$\Box$

\medskip

\noindent It follows immediately from Theorem \ref{theorem4}
that $K_{1,3}$ is the only subcubic tree $T$ with $\gamma(T)=\gamma_{e,f}^*(T)$.

Our next result gives lower bounds on the fractional exponential domination number in more general settings.
Again its proof relies on linear programming arguments.

\begin{theorem}\label{theoremlb}
Let $G$ be a graph of order $n$, maximum degree $\Delta$, and diameter $d$.
\begin{enumerate}[(i)]
\item $\gamma_{e,f}^*(G)\geq \frac{d+3}{6}$.
\item If $\Delta=3$, then $\gamma_{e,f}^*(G)\geq \frac{n}{2+3d}$.
\item If $\Delta\geq 4$, then $\gamma_{e,f}^*(G)\geq
\left(\frac{\left(\frac{\Delta-1}{2}\right)-1}{\Delta\left(\frac{\Delta-1}{2}\right)^d-3}\right)n$.
\end{enumerate}
\end{theorem}
{\it Proof:} (i) Let $P:x_0\ldots x_d$ be a shortest path of length $d$ in $G$.
Let $(y(u))_{u\in V(G)}$ be such that
$$y(u)
=
\begin{cases}
\frac{1}{3} & \mbox{, if $u$ is an endvertex of $P$,}\\
\frac{1}{6} & \mbox{, if $u$ is an internal vertex of $P$, and}\\
0 & \mbox{, if $u$ is in $V(G)\setminus V(P)$}.
\end{cases}
$$
Let $u\in V(G)$.

If $u=x_i$ for some $i\in [d-1]$, then
\begin{eqnarray*}
&& \sum\limits_{v\in V(G)}\left(\frac{1}{2}\right)^{{\rm dist}_{G}(u,v)-1}\cdot y(v)\\
&=& \sum\limits_{v\in V(P)}\left(\frac{1}{2}\right)^{{\rm dist}_{P}(u,v)-1}\cdot y(v)\\
&=&
\underbrace{
\left(\frac{1}{2}\right)^{i-1}\cdot \frac{1}{3}
+
\sum_{j=1}^{i-1}
\left(\frac{1}{2}\right)^{i-j-1}\cdot \frac{1}{6}}_{=\frac{1}{3}}
+
\underbrace{\left(\frac{1}{2}\right)^{-1}\cdot \frac{1}{6}}_{=\frac{1}{3}}
+
\underbrace{\sum_{j=1}^{d-i-1}
\left(\frac{1}{2}\right)^{j-1}\cdot \frac{1}{6}
+
\left(\frac{1}{2}\right)^{d-i-1}\cdot \frac{1}{3}}_{=\frac{1}{3}}\\
&=&1.
\end{eqnarray*}
Similarly, if $i\in \{ 0,d\}$, then we obtain
$\sum\limits_{v\in V(G)}\left(\frac{1}{2}\right)^{{\rm dist}_{G}(u,v)-1}\cdot y(v)=1$.

Now, let $u\in V(G)\setminus V(P)$.
By Lemma 4 in \cite{ddems},
there is a vertex $u'\in V(P)$ with
${\rm dist}_G(u,v)\geq {\rm dist}_G(u',v)$
for every vertex $v$ of $P$.
This implies
\begin{eqnarray*}
 \sum\limits_{v\in V(G)}\left(\frac{1}{2}\right)^{{\rm dist}_G(u,v)-1}\cdot y(v)
&= & \sum\limits_{v\in V(P)}\left(\frac{1}{2}\right)^{{\rm dist}_G(u,v)-1}\cdot y(v)\\
&\leq & \sum\limits_{v\in V(P)}\left(\frac{1}{2}\right)^{{\rm dist}_G(u',v)-1}\cdot y(v)\\
&\leq & 1.
\end{eqnarray*}
Altogther, we obtain that $(y(u))_{u\in V(G)}$
is a feasible solution for the dual of the linear program (\ref{ep})
(cf. (\ref{ed}) with ``$T$'' replaced by ``$G$''), and, by weak duality,
$\gamma_{e,f}^*(G)\geq \sum\limits_{v\in V(G)}y(v)=\frac{d+3}{6}$.

\medskip

\noindent (ii) and (iii) Let $(y(u))_{u\in V(G)}$ be such that $y(u)=y$ for every vertex $u$ of $G$ and some $y>0$.
Let $u\in V(G)$.
Since there are at most $\Delta(\Delta-1)^{i-1}$ vertices at distance $i$ from $u$ for $i\in [d]$, we obtain
\begin{eqnarray*}
\sum\limits_{v\in V(G)}\left(\frac{1}{2}\right)^{{\rm dist}_G(u,v)-1}\cdot y(v)
& \leq & 2y+\sum_{i=1}^d\Delta(\Delta-1)^{i-1}\left(\frac{1}{2}\right)^{i-1}y\\
& = &
\begin{cases}
(2+3d)y & \mbox{, if $\Delta=3$, and}\\
\frac{\Delta\left(\frac{\Delta-1}{2}\right)^d-3}{\left(\frac{\Delta-1}{2}\right)-1}y & \mbox{, if $\Delta\geq 4$.}
\end{cases}
\end{eqnarray*}
If $\Delta=3$, then choosing $y=\frac{1}{2+3d}$ yields a feasible solution of the dual of (\ref{ep}),
which implies $\gamma_{e,f}^*(G)\geq \frac{n}{2+3d}$.
Similarly, if $\Delta\geq 4$, then choosing $y=
\frac{\left(\frac{\Delta-1}{2}\right)-1}{\Delta\left(\frac{\Delta-1}{2}\right)^d-3}$ yields
$\gamma_{e,f}^*(G)\geq \left(\frac{\left(\frac{\Delta-1}{2}\right)-1}{\Delta\left(\frac{\Delta-1}{2}\right)^d-3}\right)n$.
$\Box$

\medskip

\noindent Since $\max\left\{ \frac{d+3}{6},\frac{n}{2+3d}\right\}\geq \frac{1}{6}\left(\sqrt{2n+\frac{49}{36}}+\frac{7}{6}\right)$, the following corollary is immediate.

\begin{corollary}\label{corollary2}
If $G$ is a subcubic graph of order $n$, then
$\gamma_{e,f}^*(G)\geq \frac{1}{6}\left(\sqrt{2n+\frac{49}{36}}+\frac{7}{6}\right)$.
\end{corollary}
Bessy et al.~\cite{bor2} proved that
$\gamma_e(G)\geq \frac{n(G)}{6\log_2(n(G)+2)+4}$
for a connected cubic graph $G$,
which is best possible up to the exponent of the $\log$-term in the denominator.
It is conceivable that similar lower bounds hold
for $\gamma_e^*(G)$ or even $\gamma_{e,f}^*(G)$,
which would greatly improve Corollary \ref{corollary2}.

\medskip

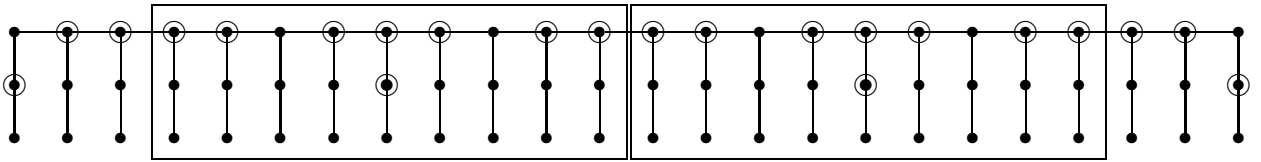
\begin{figure}[H]
\begin{center}
\unitlength 0.7mm 
\linethickness{0.4pt}
\ifx\plotpoint\undefined\newsavebox{\plotpoint}\fi 
\begin{picture}(237,30)(0,0)
\put(125,5){\circle*{2}}
\put(35,5){\circle*{2}}
\put(165,5){\circle*{2}}
\put(75,5){\circle*{2}}
\put(205,5){\circle*{2}}
\put(115,5){\circle*{2}}
\put(5,5){\circle*{2}}
\put(135,5){\circle*{2}}
\put(45,5){\circle*{2}}
\put(175,5){\circle*{2}}
\put(85,5){\circle*{2}}
\put(215,5){\circle*{2}}
\put(15,5){\circle*{2}}
\put(145,5){\circle*{2}}
\put(55,5){\circle*{2}}
\put(185,5){\circle*{2}}
\put(95,5){\circle*{2}}
\put(225,5){\circle*{2}}
\put(25,5){\circle*{2}}
\put(155,5){\circle*{2}}
\put(65,5){\circle*{2}}
\put(195,5){\circle*{2}}
\put(105,5){\circle*{2}}
\put(235,5){\circle*{2}}
\put(125,15){\circle*{2}}
\put(35,15){\circle*{2}}
\put(165,15){\circle*{2}}
\put(75,15){\circle*{2}}
\put(205,15){\circle*{2}}
\put(115,15){\circle*{2}}
\put(5,15){\circle*{2}}
\put(135,15){\circle*{2}}
\put(45,15){\circle*{2}}
\put(175,15){\circle*{2}}
\put(85,15){\circle*{2}}
\put(215,15){\circle*{2}}
\put(15,15){\circle*{2}}
\put(145,15){\circle*{2}}
\put(55,15){\circle*{2}}
\put(185,15){\circle*{2}}
\put(95,15){\circle*{2}}
\put(225,15){\circle*{2}}
\put(25,15){\circle*{2}}
\put(155,15){\circle*{2}}
\put(65,15){\circle*{2}}
\put(195,15){\circle*{2}}
\put(105,15){\circle*{2}}
\put(235,15){\circle*{2}}
\put(5,25){\circle*{2}}
\put(135,25){\circle*{2}}
\put(45,25){\circle*{2}}
\put(175,25){\circle*{2}}
\put(85,25){\circle*{2}}
\put(215,25){\circle*{2}}
\put(15,25){\circle*{2}}
\put(145,25){\circle*{2}}
\put(55,25){\circle*{2}}
\put(185,25){\circle*{2}}
\put(95,25){\circle*{2}}
\put(225,25){\circle*{2}}
\put(25,25){\circle*{2}}
\put(155,25){\circle*{2}}
\put(65,25){\circle*{2}}
\put(195,25){\circle*{2}}
\put(105,25){\circle*{2}}
\put(235,25){\circle*{2}}
\put(5,25){\line(0,-1){20}}
\put(135,25){\line(0,-1){20}}
\put(45,25){\line(0,-1){20}}
\put(175,25){\line(0,-1){20}}
\put(85,25){\line(0,-1){20}}
\put(215,25){\line(0,-1){20}}
\put(15,25){\line(0,-1){20}}
\put(145,25){\line(0,-1){20}}
\put(55,25){\line(0,-1){20}}
\put(185,25){\line(0,-1){20}}
\put(95,25){\line(0,-1){20}}
\put(225,25){\line(0,-1){20}}
\put(25,25){\line(0,-1){20}}
\put(155,25){\line(0,-1){20}}
\put(65,25){\line(0,-1){20}}
\put(195,25){\line(0,-1){20}}
\put(105,25){\line(0,-1){20}}
\put(235,25){\line(0,-1){20}}
\put(125,15){\line(0,-1){10}}
\put(35,15){\line(0,-1){10}}
\put(165,15){\line(0,-1){10}}
\put(75,15){\line(0,-1){10}}
\put(205,15){\line(0,-1){10}}
\put(115,15){\line(0,-1){10}}
\put(125,25){\circle*{2}}
\put(35,25){\circle*{2}}
\put(165,25){\circle*{2}}
\put(75,25){\circle*{2}}
\put(205,25){\circle*{2}}
\put(115,25){\circle*{2}}
\put(125,24){\line(0,-1){9}}
\put(35,24){\line(0,-1){9}}
\put(165,24){\line(0,-1){9}}
\put(75,24){\line(0,-1){9}}
\put(205,24){\line(0,-1){9}}
\put(115,24){\line(0,-1){9}}
\put(5,25){\line(1,0){9}}
\put(135,25){\line(1,0){9}}
\put(45,25){\line(1,0){9}}
\put(175,25){\line(1,0){9}}
\put(85,25){\line(1,0){9}}
\put(215,25){\line(1,0){9}}
\put(15,25){\line(1,0){9}}
\put(145,25){\line(1,0){9}}
\put(55,25){\line(1,0){9}}
\put(185,25){\line(1,0){9}}
\put(95,25){\line(1,0){9}}
\put(225,25){\line(1,0){9}}
\put(25,25){\line(1,0){9}}
\put(155,25){\line(1,0){9}}
\put(65,25){\line(1,0){9}}
\put(195,25){\line(1,0){9}}
\put(105,25){\line(1,0){9}}
\put(125,25){\line(1,0){10}}
\put(35,25){\line(1,0){10}}
\put(165,25){\line(1,0){10}}
\put(75,25){\line(1,0){10}}
\put(205,25){\line(1,0){10}}
\put(5,15){\circle{4}}
\put(15,25){\circle{4}}
\put(25,25){\circle{4}}
\put(125,25){\circle{4}}
\put(35,25){\circle{4}}
\put(135,25){\circle{4}}
\put(45,25){\circle{4}}
\put(155,25){\circle{4}}
\put(65,25){\circle{4}}
\put(165,25){\circle{4}}
\put(75,25){\circle{4}}
\put(165,15){\circle{2}}
\put(75,15){\circle{2}}
\put(165,15){\circle{4}}
\put(75,15){\circle{4}}
\put(175,25){\circle{4}}
\put(85,25){\circle{4}}
\put(195,25){\circle{4}}
\put(105,25){\circle{4}}
\put(205,25){\circle{4}}
\put(115,25){\circle{4}}
\put(215,25){\circle{4}}
\put(225,25){\circle{4}}
\put(235,15){\circle{4}}
\put(121,1){\framebox(89,29)[cc]{}}
\put(31,1){\framebox(89,29)[cc]{}}
\put(116,25){\line(1,0){10}}
\end{picture}
\end{center}
\caption{A family of trees $T$ for which
$\gamma_e(T)=\gamma(T)$ and
$\frac{\gamma_e^*(T)}{\gamma_e(T)}\leq \frac{8}{9}+o(n(T))$.
The circled vertices indicate some porous exponential dominating set.}\label{fign5}
\end{figure}

\noindent Our next goal is a characterization of the subcubic trees $T$ with $\gamma(T)=\gamma_e(T)$.
As we have observed in the introduction,
no efficient algorithms are known to
determine the exponential domination number of general trees
or the porous exponential number of subcubic trees.
Therefore, it seems difficult to extend our characterization to all trees,
or to characterize the subcubic trees $T$
with $\gamma(T)=\gamma_e^*(T)$.
The family of trees
obtained by repeating the pattern
indicated in Figure \ref{fign5}
shows that the subcubic trees $T$
with $\gamma(T)=\gamma_e^*(T)$
form a proper subset of those with
$\gamma(T)=\gamma_e(T)$.

\noindent Let $G$ be a graph.
For a vertex $x$ of $G$, let $\tau_G(x)$ be the minimum real value $\tau$ such that
there is a set $D$ of vertices of $G$ with
\begin{itemize}
\item $|D|<\gamma_e(G)$,
\item $x\not\in D$,
\item $w_{(G,D)}(u)+\left(\frac{1}{2}\right)^{{\rm dist}_{G-D}(x,u)}\cdot\tau\geq 1$ for every vertex $u$ in $V(G)\setminus D$.
\end{itemize}
Now, we define three operations on trees.
Let $T$ and $T'$ be two trees.
\begin{itemize}
\item {\bf Operation 1}

$T$ arises from $T'$ by {\it applying Operation 1} if $T$ has an endvertex $y$ with neighbor $x$ such that
$T'=T-y$, and
$x$ belongs to some minimum dominating set of $T'$.
\item {\bf Operation 2}

$T$ arises from $T'$ by {\it applying Operation 2} if $T$ contains a path $xyz$ such that
$\tau_{T'}(x)>1$ or $\gamma(T',V(T')\setminus \{ x\})<\gamma(T')$,
where
$y$ has degree $2$ in $T$,
$z$ is an endvertex of $T$, and
$T'=T-\{ y,z\}$.
\item {\bf Operation 3}

$T$ arises from $T'$ by {\it applying Operation 3} if $T$ contains a path $wxyz$ such that
$x$ and $y$ have degree $2$ in $T$,
$z$ is an endvertex of $T$,
$T'=T-\{ x,y,z\}$, and
$\tau_{T'}(w)>\frac{1}{2}$.
\end{itemize}
Let ${\cal T}$ be the family of subcubic trees that are obtained from $P_1$ by applying finite sequences of the above three operations.

\begin{lemma}\label{lemma2}
If $T'$ is a subtree of a subcubic tree $T$, then $\gamma_e(T')\leq \gamma_e(T)$.
\end{lemma}
{\it Proof:}
By an inductive argument, it suffices to consider the case that $T'=T-y$,
where $y$ is an endvertex of $T$.
Let $x$ be the neighbor of $y$.
Let $D$ be a minimum exponential dominating set of $T$.
If $y\not\in D$, then $D$ is also an exponential dominating set of $T'$.
If $y\in D$, then $x\not\in D$, because $D$ is minimum.
Let $D'=(D\setminus \{ y\})\cup \{ x\}$.
Suppose that there is some vertex $u$ with $w_{(T,D')}(u)<1$.
Clearly, $u\not=x$.
Let $e$ be the edge of the path $P$ between $u$ and $x$ that is incident with $x$.
Let $T_x$ be the component of $T-e$ that contains $x$, and let $D_x=D\cap V(T_x)$.
Since $w_{(T,D')}(u)<w_{(T,D)}(u)$, the path $P$ does not intersect $D$.
This implies
\begin{eqnarray*}
w_{(T,D')}(u)
& = &
w_{(T,D)}(u)-w_{(T,D_x)}(u)+w_{(T,\{ x\})}(u)\\
& = &
w_{(T,D)}(u)-\left(\frac{1}{2}\right)^{{\rm dist}_T(u,x)}\cdot w_{(T_x,D_x)}(x)+\left(\frac{1}{2}\right)^{{\rm dist}_T(u,x)-1}.
\end{eqnarray*}
Since $w_{(T,D')}(u)<w_{(T,D)}(u)$, this implies $w_{(T_x,D_x)}(x)>2$,
which contradicts Lemma \ref{lemma1}.
Hence, $D'$ is an exponential dominating set of $T'$.
Altogether, we obtain $\gamma_e(T')\leq \gamma_e(T)$.
$\Box$

\begin{lemma}\label{lemma3}
If $T\in {\cal T}$, then $\gamma_e(T)=\gamma(T)$.
\end{lemma}
{\it Proof:}
Note that $\gamma(P_1)=\gamma_e(P_1)$.
By an inductive argument, it suffices to show that $\gamma(T)=\gamma_e(T)$
for every tree $T$ that arises from some tree $T'$ with $\gamma(T')=\gamma_e(T')$
by applying one of the above three operations.

First, let $T$ arise from $T'$ by applying Operation 1.
Since $x$ belongs to some minimum dominating set of $T'$, we have $\gamma(T')=\gamma(T)$.
By Lemma \ref{lemma2}, we obtain
$$\gamma(T)=\gamma(T')=\gamma_e(T')\leq \gamma_e(T)\leq \gamma(T),$$
which implies $\gamma_e(T)=\gamma(T)$.

Next,
let $T$ arise from $T'$ by applying Operation 2.

First, we assume that $\tau_{T'}(x)>1$.
By Lemma \ref{lemma2}, we have $\gamma_e(T')\leq \gamma_e(T)$.
Suppose that $\gamma_e(T')=\gamma_e(T)$.
Let $D$ be a minimum exponential dominating set of $T$.
By Lemma \ref{lemma1}, the set $D$ must contain either $y$ or $z$.
Clearly, we may assume $y\in D$ and $z\not\in D$.
Let $D'=D\setminus \{ y\}$.
Since $|D'|<\gamma_e(T')$, the set $D'$ is not an exponential dominating set of $T'$,
which implies that $x\not\in D'$.
Since
\begin{eqnarray*}
w_{(T,D)}(u)
& = & w_{(T,D')}(u)+\left(\frac{1}{2}\right)^{{\rm dist}_{T-D'}(u,y)-1}\\
& = & w_{(T,D')}(u)+\left(\frac{1}{2}\right)^{{\rm dist}_{T'-D'}(u,x)}\cdot 1\\
& \geq & 1,
\end{eqnarray*}
for every vertex $u\in V(T')\setminus D'$,
we obtain the contradiction that $\tau_{T'}(x)\leq 1$.
Hence, $\gamma_e(T')+1\leq\gamma_e(T)$.
Note that $\gamma(T)\leq \gamma(T')+1$.
Now,
$$\gamma(T)\leq \gamma(T')+1=\gamma_e(T')+1\leq \gamma_e(T)\leq \gamma(T),$$
which implies $\gamma_e(T)=\gamma(T)$.

Next, we assume that $\gamma(T',V(T')\setminus \{ x\})<\gamma(T')$.
Let $D'$ be a set of vertices of $T'$ with $|D'|=\gamma(T',V(T')\setminus \{ x\})$
such that every vertex in $(V(T')\setminus \{ x\})\setminus D'$ has a neighbor in $D'$.
Since $D'\cup \{ y\}$ is a dominating set of $T$, we obtain $\gamma(T)\leq  \gamma(T',V(T')\setminus \{ x\})+1\leq \gamma(T')$.
By Lemma \ref{lemma2}, we obtain
$$\gamma_e(T)\leq \gamma(T)\leq \gamma(T')=\gamma_e(T')\leq \gamma_e(T),$$
which implies $\gamma_e(T)=\gamma(T)$.

Next,
let $T$ arise from $T'$ by applying Operation 3.
Clearly, $\gamma(T)=\gamma(T')+1$.
Suppose that $\gamma_e(T)\leq \gamma_e(T')$.
Let $D$ be a minimum exponential dominating set of $T$.
By Lemma \ref{lemma1}, the set $D$ must contain either $y$ or $z$.
Clearly, we may assume $y\in D$ and $z\not\in D$.
Arguing similarly as in the proof of Lemma \ref{lemma2},
we may assume that $x\not\in D$.
Let $D'=D\setminus \{ y\}$.
Since $|D'|<\gamma_e(T')$, the set $D'$ is not an exponential dominating set of $T'$,
which implies that $w\not\in D'$.
Since
\begin{eqnarray*}
w_{(T,D)}(u)
& = & w_{(T,D')}(u)+\left(\frac{1}{2}\right)^{{\rm dist}_{T-D'}(u,y)-1}\\
& = & w_{(T,D')}(u)+\left(\frac{1}{2}\right)^{{\rm dist}_{T'-D'}(u,w)}\cdot \frac{1}{2}\\
& \geq & 1,
\end{eqnarray*}
for every vertex $u\in V(T')\setminus D'$,
we obtain the contradiction that $\tau_{T'}(w)\leq \frac{1}{2}$.
Hence, $\gamma_e(T')+1\leq\gamma_e(T)$.
Now,
$$\gamma(T)=\gamma(T')+1=\gamma_e(T')+1\leq\gamma_e(T)\leq\gamma(T),$$
which implies $\gamma_e(T)=\gamma(T)$.
$\Box$

\begin{theorem}\label{theorem1}
If $T$ is a subcubic tree, then $\gamma(T)=\gamma_e(T)$ if and only if $T\in {\cal T}$.
\end{theorem}
{\it Proof:} Lemma \ref{lemma3} implies the sufficiency.
In order to prove the necessity,
suppose that $T$ is a subcubic tree of minimum order
such that $\gamma(T)=\gamma_e(T)$ and $T\not\in {\cal T}$.
Considering three applications of Operation 1 to $P_1$ implies $P_2,P_3,K_{1,3}\in {\cal T}$,
that is, ${\cal T}$ contains all subcubic trees of diameter at most $2$.
Hence, the tree $T$ has diameter at least $3$.
Let $v$ be an endvertex of a longest path in $T$.
The vertex $v$ has a unique neighbor $u$ in $T$.

\begin{claim}\label{claim1}
The vertex $u$ has degree $2$.
\end{claim}
{\it Proof of Claim \ref{claim1}:}
Suppose that $u$ has degree $3$ in $T$.
This implies that $u$ has a neighbor $w$ that is an endvertex of $T$ distinct from $v$.
Let $T'=T-w$.
Clearly, $u$ belongs to some minimum dominating set of $T'$,
which implies $\gamma(T)=\gamma(T')$.
Arguing as above, we obtain that $T'$ has a minimum exponential dominating set $D'$ that contains $u$.
Since $D'$ is also an exponential dominating set of $T$, we obtain $\gamma_e(T)\leq \gamma_e(T')$.
By Lemma \ref{lemma2}, we have $\gamma_e(T)=\gamma_e(T')$.
Now, $\gamma(T')=\gamma(T)=\gamma_e(T)=\gamma_e(T')$.
By the choice of $T$, we obtain $T'\in {\cal T}$.
Since $T$ arises from $T'$ by applying Operation 1,
we obtain $T\in {\cal T}$,
which is a contradiction. $\Box$

\medskip

\noindent Let $x$ be the neighbor of $u$ distinct from $v$.
Let $T''=T-\{ u,v\}$.

\begin{claim}\label{claim2}
$\tau_{T''}(x)>1$ or $\gamma(T'',V(T'')\setminus \{ x\})<\gamma(T'')$.
\end{claim}
{\it Proof of Claim \ref{claim2}:}
Suppose that $\tau_{T''}(x)\leq 1$ and $\gamma(T'',V(T'')\setminus \{ x\})\geq \gamma(T'')$.
Arguing as above, we obtain that the first condition implies $\gamma_e(T)=\gamma_e(T'')$,
and that the second condition implies $\gamma(T)=\gamma(T'')+1$.
Now, we obtain the contradiction $\gamma_e(T)=\gamma_e(T'')\leq \gamma(T'')<\gamma(T)$.
$\Box$

\medskip

\noindent If $\gamma_e(T'')=\gamma(T'')$, then, by the choice of $T$, we have $T''\in {\cal T}$,
and, by Claim \ref{claim2}, the tree $T$ arises from $T''$ by applying Operation 2,
which implies the contradiction $T\in {\cal T}$.
Hence, we may assume that $\gamma_e(T'')<\gamma(T'')$.

\medskip

\noindent
According to Claim \ref{claim2}, we consider two cases.

\medskip

\noindent {\bf Case 1} {$\tau_{T''}(x)>1$.}

\medskip

\noindent Clearly, $\gamma_e(T)\leq \gamma_e(T'')+1$, and we obtain
$$\gamma(T)=\gamma_e(T)\leq \gamma_e(T'')+1\leq \gamma(T'')\leq \gamma(T),$$
which implies $\gamma(T'')=\gamma(T)=\gamma_e(T'')+1$.
If $x$ has degree $3$, then, by the choice of $v$,
either $x$ has a neighbor that is an endvertex
or there is a path $v'u'x$, where $v'$ is an endvertex that is distinct from $v$.
In both cases, $T$ has a minimum dominating set that contains $u$ and either $x$ or $u'$,
which implies the contradiction $\gamma(T'')<\gamma(T)$.
Hence, the vertex $x$ has degree $2$.
Let $y$ be the neighbor of $x$ that is distinct from $u$.
Let $T'''=T-\{ x,u,v\}$.
Suppose that $\tau_{T'''}(y)\leq \frac{1}{2}$.
Arguing as above, this implies $\gamma_e(T)=\gamma_e(T''')$.
By Lemma \ref{lemma2}, we obtain the contradiction
$\gamma_e(T''')=\gamma_e(T)=\gamma_e(T'')+1\geq \gamma_e(T''')+1$.
Hence, $\tau_{T'''}(y)>\frac{1}{2}$,
which implies $\gamma_e(T)=\gamma_e(T''')+1$.
Since $\gamma(T)=\gamma(T''')+1$,
we obtain $\gamma_e(T''')=\gamma_e(T)-1=\gamma(T)-1=\gamma(T''')$.
By the choice of $T$, this implies that $T'''\in {\cal T}$,
and that the tree $T$ arises from $T'''$ by applying Operation 3,
which implies the contradiction $T\in {\cal T}$.

\medskip

\noindent {\bf Case 2} {$\gamma(T'',V(T'')\setminus \{ x\})<\gamma(T'')$.}

\medskip

\noindent Clearly, $\gamma(T'')=\gamma(T)$.
As in Case 1, this implies that $x$ has degree $2$.
Let $y$ and $T'''$ be as in Case 1.
Suppose that $\tau_{T'''}(y)\leq \frac{1}{2}$.
Similarly as above, this implies the contradiction
$\gamma(T''')\geq \gamma_e(T''')\geq\gamma_e(T)=\gamma(T)=\gamma(T''')+1$.
Hence, $\tau_{T'''}(y)>\frac{1}{2}$,
which implies $\gamma_e(T)=\gamma_e(T''')+1$.
Again $\gamma(T)=\gamma(T''')+1$, and
we obtain $\gamma_e(T''')=\gamma_e(T)-1=\gamma(T)-1=\gamma(T''')$.
By the choice of $T$, this implies that $T'''\in {\cal T}$,
and that the tree $T$ arises from $T'''$ by applying Operation 3,
which implies the contradiction $T\in {\cal T}$.

This completes the proof. $\Box$

\medskip

\noindent A drawback of the above characterization is the use of the values $\tau_G(u)$ and conditions such as
``$\gamma(T',V(T')\setminus \{ x\}))<\gamma(T')$''
in the definition of Operation 2.
It is conceivable that these technical complications can be eliminated, and that a completely explicit (constructive) characterization is possible.

\end{document}